\documentclass[11pt]{article}
   \usepackage{amsfonts}
\usepackage{latexsym,tikz,mathdots}
\usepackage{amsmath}
 
\title{Multiple Dedekind Zeta Functions}
\author{Ivan Emilov Horozov}
\date{February 2, 2008}

\newcommand \nc {\newcommand}
\nc \proof {\noindent {\em{Proof.\/ }}} \nc \qed {$\Box$\hfill}
\newtheorem{theorem}{Theorem}[section]
\newtheorem{lemma}[theorem]{Lemma}
\newtheorem{proposition}[theorem]{Proposition}
\newtheorem{corollary}[theorem]{Corollary}
\newtheorem{definition}[theorem]{Definition}
\newtheorem{example}[theorem]{Example}
\newtheorem{remark}[theorem]{Remark}
\newtheorem{conjecture}[theorem]{Conjecture}
\newtheorem{question}[theorem]{Question}
\nc \bth[1] {\begin{theorem}\label{t#1} } \nc \ble[1]
{\begin{lemma}\label{l#1} } \nc \bpr[1]
{\begin{proposition}\label{p#1} } \nc \bco[1]
{\begin{corollary}\label{c#1} } \nc \bde[1]
{\begin{definition}\label{d#1}\rm } \nc \bex[1]
{\begin{example}\label{e#1}\rm } \nc \bre[1]
{\begin{remark}\label{r#1}\rm } \nc \bcon[1]
{\begin{conjecture}\label{con#1}\rm } \nc \bque[1]
{\begin{question}\label{que#1}\rm }
\nc {\eth} { \end{theorem} } \nc {\ele} { \end{lemma} } \nc
{\epr}{ \end{proposition} } \nc {\eco} { \end{corollary} } \nc
{\ede} {\end{definition} } \nc {\eex} { \end{example} } \nc {\ere}
{\end{remark} } \nc {\econ} { \end{conjecture} } \nc {\eque}
{\end{question} }

\def \Q {{\mathbb Q}}

\begin{document}

\title{{\LARGE\bf{An Exponential Diophantine Equation -\\ One Order Higher than Fermat's Equation}}}

\author{
I. ~Horozov
\thanks{E-mail: horozov@bcc.cuny.edu}
\\ \hfill\\ \normalsize \textit{Department of Mathematics and Computer Science,}\\
\normalsize \textit{City University of New York, Bronx Community College}\\
\normalsize \textit{2155 University Avenue, Bronx,
New York 10453 , U.S.A.}\\
}
\date{}
\maketitle

{\it{To {\bf{Yuri Manin}} with respect and admiration.}}

\vspace{.2cm}
\begin{abstract}
We formulate an exponential Diophantine equation, which is is some sense one order higher that Fermat's Last Theorem. 
We also give three examples of solutions to this exponential Diophantine equation and formulate a conjecture.
\end{abstract}




\
\section{An Exponential Diophantine Equation}
Fermat Last Theorem was an open conjecture for about 300 years. It was proven not long ago by Andrew Wiles \cite{W}.
\begin{theorem} (Fermat's Last Theorem; Wiles)
For an integer $n>2$ there are no nontrivial integer  solutions of
\[x^n+y^n=z^n.\] 
Where by trivial we mean one of the variables $x,y,z$ to be zero.
\end{theorem}
Earlier proofs, for a particular  exponent $n$ or for a large class of exponents $n$, were based on a algebraic number theory. One can look at the book by Ireland and Rosen \cite{IR} for Euler's proof for $n=3,4$, and \cite{K} for Kummer's proof for a large class of exponents $n$.
Wiles proof was based on proving Tanyiama-Shimura conjecture that all elliptic curves over $\Q$ are modular.

More recently, a slight modification of the formulation of the Fermat Last Theorem was related to Hilbert modular surfaces \cite{Mladen}.
One may consider the following Diophantine equation:
\begin{equation}
\label{eq2}
x^k+y^m=z^n,
\end{equation}
where $k$, $m$ and $n$ could be different positive integers, (see \cite{Kr}).

We generalize Equation \eqref{eq2} so that the operation are in some sense one order of complexity higher. Instead of addition we consider multiplication and instead of multiplication we consider exponentiation. Note that for a positive real number $\alpha$, we can define $\alpha^\alpha$. However, this operation is not associative. For example, in general
\[\alpha^{\left(\alpha^\alpha\right)}
\neq
\left(\alpha^\alpha\right)^\alpha.
\]
By successive exponentiation we mean 
\[\alpha^{\alpha^\alpha}
=
\alpha^{\left(\alpha^\alpha\right)}
\]
and
\[\alpha^{\alpha^{\alpha^\alpha}}
=
\alpha^{\left(\alpha^{\left(\alpha^\alpha\right)}\right)}
\]

\begin{conjecture}
Let $\alpha$, $\beta$ and $\gamma$ are positive real algebraic numbers. 
Consider an exponential Diophantine equation 
\begin{equation}
\label{eq exp}
\alpha^{\iddots^\alpha}\times\, \beta^{\iddots^\beta} = \gamma^{\iddots^{\gamma}},
\end{equation}
where $\alpha$ is repeated $k$-times, $\beta$ is repeated $m$-times, and $\gamma$ is repeated $n$-times, for $k$, $m$ and $n$ greater that $1$.
Then, the Equation \eqref{eq exp} has only finitely many solutions for any fixed $k$, $m$ and $n$.
\end{conjecture}
A number $\alpha$ is called {\it{algebraic}} if there is a polynomial $f(x)$ with rational coefficients such that $f(\alpha)=0$.
For an introduction to this topic one could consider the book by Ireland and Rosen \cite{IR} and Lang \cite{L}.
The trivial solutions are: $\alpha=1$ and $\beta=\gamma$, or $\beta=1$ and $\alpha=\gamma$.

\section{Examples}
We present several examples of solutions to this exponential Diophantine equation.
\begin{theorem}
(a)\[2^{2^2}\times 2^{2^2}=4^4;\]
(b) \[\left(\frac{1}{2}\right)^{\left(\frac{1}{2}\right)^{\left(\frac{1}{2}\right)}}
\times
 \left(\frac{1}{2}\right)^{\left(\frac{1}{2}\right)^{\left(\frac{1}{2}\right)}}
 =
 \left(\frac{1}{4}\right)^{\left(\frac{1}{4}\right)^{\left(\frac{1}{4}\right)}};\]
 (c) 
  \[
\left(\frac{1}{\sqrt{2}}\right)^{\left(\frac{1}{\sqrt{2}}\right)}
\times
\left(\frac{1}{\sqrt{2}}\right)^{\left(\frac{1}{\sqrt{2}}\right)}
=
 \left(\frac{1}{2}\right)^{\left(\frac{1}{2}\right)^{\left(\frac{1}{2}\right)}}.\]
\end{theorem}
\proof 
For part (a), we have
\[2^{2^2}\times 2^{2^2}=2^4\times 2^4=2^8=2^{2\times 4}=\left(2^2\right)^4=4^4.\]
We compute the left hand side of part (b).
\[\left(\frac{1}{2}\right)^{\left(\frac{1}{2}\right)}=\frac{1}{\sqrt{2}}.\]
Then
\begin{equation}
\label{eq1}
  \left(\frac{1}{2}\right)^{\left(\frac{1}{2}\right)^{\left(\frac{1}{2}\right)}}=2^{\left(-\frac{1}{\sqrt{2}}\right)}
\end{equation}
And the left hand side of part (b) becomes
  \[
  \left(\frac{1}{2}\right)^{\left(\frac{1}{2}\right)^{\left(\frac{1}{2}\right)}}
\times
 \left(\frac{1}{2}\right)^{\left(\frac{1}{2}\right)^{\left(\frac{1}{2}\right)}}
=
\left( 
\left(\frac{1}{2}\right)^{\left(\frac{1}{2}\right)^{\left(\frac{1}{2}\right)}}
\right)^2
=
2^{-\frac{2}{\sqrt{2}}}
=
2^{-\sqrt{2}}.\]
For the right hand side of part (b), we have
\[
\left(\frac{1}{4}\right)^{\left(\frac{1}{4}\right)}
=
4^{-\left(\frac{1}{4}\right)}
=
2^{-2\left(\frac{1}{4}\right)}
=
2^{-{\left(\frac{1}{2}\right)}}=\frac{1}{\sqrt{2}}
\]
Then
\[
\left(\frac{1}{4}\right)^{\left(\frac{1}{4}\right)^{\left(\frac{1}{4}\right)}}
=
\left(\frac{1}{4}\right)^{\frac{1}{\sqrt{2}}}
=
4^{-\frac{1}{\sqrt{2}}}
=
2^{-2{\frac{1}{\sqrt{2}}}}
=
2^{-\sqrt{2}}.
\]

For part (c) we have
\[
\left(\frac{1}{\sqrt{2}}\right)^{\left(\frac{1}{\sqrt{2}}\right)}
=\sqrt{2}^{\left(-\frac{1}{\sqrt{2}}\right)}
=
2^{\left(-\frac{1}{2}\cdot\frac{1}{\sqrt{2}}\right)}
\]
Then
\[  
\left(\frac{1}{\sqrt{2}}\right)^{\left(\frac{1}{\sqrt{2}}\right)}
\times
\left(\frac{1}{\sqrt{2}}\right)^{\left(\frac{1}{\sqrt{2}}\right)}
=
\left(
2^{\left(-\frac{1}{2}\cdot\frac{1}{\sqrt{2}}\right)}
\right)^2
=
2^{\left(-\frac{1}{\sqrt{2}}\right)}
\]
For the right hand side of part (c), we have
\[  \left(\frac{1}{2}\right)^{\left(\frac{1}{2}\right)^{\left(\frac{1}{2}\right)}}=2^{\left(-\frac{1}{\sqrt{2}}\right)}\]
from Equation \ref{eq1}, which agrees with the left hand side of part (c).
\qed

\section*{Acknowledgements:} I would like to thank Alexander Kosyak for suggestions to improve the formulation of Conjecture 1.2.
I would like also to thank Xuanyu Pan for checking the correctness of the computations and for his enthusiasm about the paper.

\renewcommand{\em}{\textrm}
\begin{small}

\end{small}


\begin{thebibliography}{BHY1}

\bibitem[BBDDDV]{Mladen}
L. Berger, G. B\"ockle, L. Demb\'el\'e, M. Dimitrov, T. Dokchitser, J. Voight:
{\em{Elliptic Curves, Hilbert Modular Forms and Galois Deformations}, Bitkh\"auser, Springer, Basel 2013, 249pp.}










\bibitem[IR]{IR}
Ireland K., Rosen, M.:
{\em{A classical introduction to modern number theory. Section edition.}
Graduate texts in Mathematics, 84 Springer-Verlag, New York, 1990, xiv+339pp ISBN:0-387-97329-X}

\bibitem[K]{K}
Kummer, E.:
{\em{Beweis des Fermat'schen Satzes der Unm\"oglichkeiten von 
$x^\lambda+y^\lambda=z^\lambda$ 
f\"ur eine enedliche Anzahl Primzahlen 
$\lambda$,}
Ernst Eduard Kummer, Collected Papers, Volume 1,
 Ed. Andr\'e Weil, Springer-Verlag, Berlin Heidelberg New York 1975, ISBN 3-540-06835-X}

\bibitem[Kr]{Kr}
Kraus, A.:
{\em{On the Equation
$x^p+y^q=z^r$}, A Survey, The Ramanujan Journal, 3 (1999), 315-333.}

\bibitem[L]{L}
Lang S.: {\em{Algebraic Number Theory,} 2nd ed. New York, Springer-Verlag, 1994.}
 

\bibitem[W]{W}
Wiles, A.:
{\em{Modular elliptic curves and Fermat's last theorem.} 
Ann. of Math. (2) 141 (1995), no. 3, 443-551.}

    \end{thebibliography}
\end{document}